\documentclass[12pt, a4paper]{article}
\usepackage{times}
\usepackage{color}

\newtheorem{theorem}{Theorem}[section]

\def\qed{\nopagebreak\hfill{\rule{4pt}{7pt}}}

\begin{document}

%%% \noindent {\small Date: 09-26-2023}

%%%  \noindent {\sf    \underline{Draft, 09-26-2023}}

\begin{center}

{\Large\bf   Breaking Cycles, the Odd Versus the Even
}

William Y.C. Chen

\vskip 3mm

Center for Applied Mathematics\\
Tianjin University\\
Tianjin 300072, P.R. China

\vskip 1mm

Email: chenyc@tju.edu.cn

\end{center}

\vskip 6mm

\begin{abstract} In an award-winning  expository article,
V. Pozdnyakov and J.M. Steele gave a beautiful demonstration of
the ramifications of
a basic bijection for permutations.
The aim of this note is to connect this correspondence to a
seemingly unrelated problem concerning odd cycles and even
cycles, arising in the combinatorial study of the Cayley continuants
by E. Munarini and D. Torri. In  extreme cases, one encounters
two special classes of permutations of $2n$ elements with the same
cardinality. A bijection of this appealing relation
has been found by E. Sayag. A combinatorial
study of permutations with only odd cycles has been
carried out by M. B\'ona, A. Mclennan and D. White.
We find an intermediate structure which leads to a
linkage between these two antipodal structures.
A recursive setting reveals that everything boils down
to only one trick -- breaking the cycles.
\end{abstract}



\section{Introduction}

In an award-winning exposition,
 V. Pozdnyakov and J.M. Steele     \cite{PS-2016}
  elaborated on many a facet
of a basic property of the cycle representation of
permutations, viz., the number of permutations
of $[n]=\{1,2,\ldots, n\}$ $(n\geq 2)$ for which
$1$ and $2$ occur in the same cycle equals the
number of permutations of $[n]$ for which $1$ and $2$ do not
occur in the same  cycle.
The heart of the plot lies in an
operation of breaking a cycle into two cycles.

More precisely, given a cycle containing both $1$ and $2$,
we can split it into two segments, one starting with
with $1$ and ending with the element preceding
$2$, whereas the other starting with $2$ and ending with
the element preceding $1$. Keep in mind that
a cycle can be expressed as a sequence starting with
the minimum element.

The objective of this note is to supplement
the showcase of Pozdnyakov-Steele with one more
story.
In a different scene,
we meet up with two classes of permutations
of $[2n]$ $(n\geq 1)$.
Let $A_{n}$ denote the set of permutations of $[n]$
consisting of odd cycles, let $B_{2n}$ denote the
set of permutations of $[2n]$ consisting of even cycles.
A bijection between $A_{2n}$ and $B_{2n}$
can be found in \cite[Section 6.2]{Bona-2005}.
Let $a_{n}=|A_{n}|$ and $b_{2n}=|B_{2n}|$.
As pointed out by Munarini and Torri \cite{MT-2005}, the
generating function of the Cayley continuants
specializes to the
generating functions for $a_{2n}$ and $b_{2n}$.
In fact, we have
\[
a_{2n}=b_{2n} = ((2n-1)!!)^2.
\]
 The sequence $\{a_n\}$ is
listed as $\#A000246$ in OEIS \cite{OEIS}, and the sequence
$\{b_{2n}\}$ is referred to  as  $\#A001818$. A further study of
the sequence $\{a_n\}$ can be found in   B\'ona-Mclennan-White
\cite{BMW-2000}.

We take a different avenue to provide
a combinatorial interpretation by
employing the Pozdnyakov-Steele bijection
with a twist of the roles of $1$ and $2$
in certain circumstances. As an intermediate step,
we establish the following
 correspondence.
 Let   $P_{2n}$ be the set of
permutations of $[2n]$ consisting of odd cycles except
that the element $1$ is in an even cycle.

\begin{theorem} \label{t-a-p}
There exists a bijection between $A_{2n}$ and $P_{2n}$.
\end{theorem}

\section{A bijection}

 Before presenting the proof, let us consider
how to apply the  map in Theorem \ref{t-a-p}  to
transform a permutation in $A_{2n}$ to a permutation   in $B_{2n}$.
Starting with a permutation in $A_{2n}$,
at the first step, we get a permutation with
$1$ appearing in an even cycle.
Iterating this procedure for the remaining odd cycles, we
are led to a permutation of even cycles. This proves that
$a_{2n}=b_{2n}$.

The following inductive proof is essentially a description of a
recursive algorithm.

\noindent{\it Inductive Proof of Theorem \ref{t-a-p}.}
For $n=1$, the required correspondence is
merely the only way to break the even cycle $(12)$ into
two odd cycles $(1)$ and $(2)$.

Assume that  $n>1$ and that there is a one-to-one correspondence
between $A_{2m}$ and $P_{2m}$ for $m < n$. We are going to
put together a bijection between $A_{2n}$ and $P_{2n}$.
To this end, we define
$P^{12}_{2n}$ to be the set of permutations in $P_{2n}$
such $1$ and $2$ belong to the same even cycle,
and denote by $P^{1-2}_{2n}$ the set of
permutations in $P_{2n}$ such that $1$ appears in an
even cycle but $2$ appears in an odd cycle.
Thus,
\[
    P_{2n} = P^{12}_{2n} \, \cup \, P^{1-2}_{2n}.
\]
For an even cycle containing both $1$ and $2$,
we may break it into two cycles with one containing $1$ and
the other containing $2$. Taking the parities into account,
we find that
\[
   P^{12}_{2n} \leftrightarrow   A^{1-2}_{2n}\, \cup \, Q^{1-2}_{2n},
\]
  where $A^{1-2}_{2n}$ is
  the set of permutations of $[2n]$   consisting of
  odd cycles such that $1$ and $2$ do not appear in
  the same cycle, and $Q^{1-2}_{2n}$ is the set of
  permutations of $[2n]$ such that
  $1$ and $2$ occur in different even cycles, whereas
  all other cycles are odd.

  Thus, it suffices to justify the following
  one-to-one correspondence
  \begin{equation} \label{pa}
      P^{1-2}_{2n} \, \cup \, Q^{1-2}_{2n} \leftrightarrow A^{12}_{2n},
  \end{equation}
  where $A^{12}_{2n}$ is
  the set of permutations of $[2n]$ consisting of
  odd cycles  such that $1$ and $2$   appear in
  the same cycle.
 By splitting a permutation in $A^{12}_{2n}$, we see that
\[
    A^{12}_{2n} =  P^{1-2}_{2n} \,  \cup \,  U^{1-2}_{2n},
\]
  where  $U^{1-2}_{2n}$ is the set of
  permutations of $[2n]$ such that
  $1$ is in an odd cycle, $2$ is in an even cycles and
  all other cycles are odd.

  In order to justify (\ref{pa}), we only need to
  establish the following
  correspondence
  \begin{equation} \label{QU}
    Q^{1-2}_{2n} \leftrightarrow U^{1-2}_{2n}.
  \end{equation}
  By exchanging the roles of $1$ and $2$,
  $U^{1-2}_{2n}$ can be identified with the
  set of permutations such that $1$ occurs in an
  even cycle and all other cycles are odd.

  Notice that the relation (\ref{QU})
  is nothing but a recursive statement of
  $A_{2n}\Leftrightarrow P_{2n}$. To be more specific,
  let $V^{1-2}_{2n}$ denote the set of permutations
  obtained from those $U^{1-2}_{2n}$ by exchanging $1$ and
  $2$. Assume that $\sigma$ is a permutation in $V^{1-2}_{2n}$
  and
  $C$ is the even cycle of $\sigma$
  containing $1$.

Invoking  the induction hypothesis with respect to
  all the odd cycles in $\sigma$, we
get an even cycle containing
  $2$ along with all other odd cycles, which is precisely
  a permutation in $Q^{1-2}_{2n}$.
This completes the proof.  \qed

\vskip 6mm \noindent{\large\bf Acknowledgments.}
I am grateful to Mikl\'os B\'ona, Sam Hopkins and Michael Wallner for
their valuable comments. This work was supported by the
National Science Foundation of China.

\end{document}